\documentclass[conference,twoside]{IEEEtran}

\usepackage{cite}
\usepackage{graphicx}
\usepackage{subfigure}
\usepackage{url}
\usepackage{amsmath,amssymb}
\usepackage{array}
\usepackage{eso-pic} % for watermark

\makeatletter % for watermark
\newcommand\BackgroundPicture[2]{%
  \setlength{\unitlength}{1pt}%
  default \put(0,\strip@pt\paperheight){%
  \parbox[t][\paperheight]{\paperwidth}{%
    \vfill
    \centering\includegraphics[angle=#2]{#1}
    \vfill
}}} %
\makeatother

\begin{document}

\AddToShipoutPicture{\BackgroundPicture{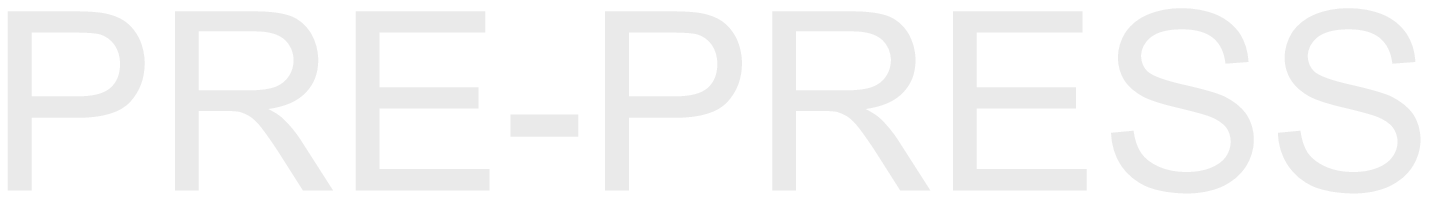}{45}} % for watermark

\newcommand{\bs}[1]{\boldsymbol{#1}}

\title{Untangling~the~SVD's~of Random~Matrix~Sample~Paths}

\author{\authorblockN{David~W.~Browne}
\authorblockA{Department of Electrical Engineering\\
University of California, Los Angeles\\
Los Angeles, California 90095-1594\\
Email: decibel@ucla.edu}
\and
\authorblockN{Michael W. Browne}
\authorblockA{Dept. of Psychology, Dept. of Statistics\\
Ohio State University\\
Columbus, Ohio 43210\\
Email: browne.4@osu.edu}
\and
\authorblockN{Michael P. Fitz}
\authorblockA{Northrop Grumman\\
Los Angeles, California\\
Email: Michael.Fitz@ngc.com}}

\maketitle
\begin{abstract}
Singular Value Decomposition (SVD) is a powerful tool for multivariate analysis. However, independent
computation of the SVD for each sample taken from a bandlimited matrix random process will result in singular value sample paths whose tangled evolution is not consistent with the structure of the underlying random process. 
The solution to this problem is developed as follows: (i) a SVD with relaxed identification conditions is proposed, (ii) an approach is formulated for computing the SVD's of two adjacent matrices in the sample path with the objective of maximizing the correlation between corresponding singular vectors of the two matrices, and (iii) an efficient algorithm is given for untangling the singular value sample paths. The algorithm gives a unique solution conditioned on the seed matrix's SVD.
Its effectiveness is demonstrated on bandlimited Gaussian random-matrix sample paths. Results are shown to be consistent with those predicted by random-matrix theory. A primary application of the algorithm is in multiple-antenna radio systems. The benefit promised by using SVD untangling in these systems is that the fading rate of the channel's SVD factors is greatly reduced so that the performance of channel estimation, channel feedback and channel prediction can be increased.
\end{abstract}

\section{Introduction}
Many observable processes in the natural world give rise to matrix random process. Examples can be found in multiple antenna radio communications \cite{Telatar1995} \cite{Foschini1998}, gene expression \cite{Alter2000} \cite{Boscolo2005}, and kinematics \cite{Maciejewski1989} among others. If these random process are bandlimited, then they may be sampled at uniform intervals and represented as discrete sample paths without loss of information. The Singular Value Decomposition (SVD) is a powerful tool for analyzing the multidimensional observations from these matrix random processes. However, computation of the SVD for each matrix-valued sample of a bandlimited random process results in singular value sample paths whose evolution is not consistent with the structure of the underlying random process. This is because the SVD computation is applied independently from sample to sample, during which strict phase and ordering identification conditions are imposed on the singular values.

This work addresses the problem of how to compute the SVD's of a random matrix sample path in a way that preserves the covariance of the underlying random process. The solution to this problem is developed as follows: (i) a SVD with relaxed identification conditions is proposed, (ii) an approach is formulated for computing the SVD's of two adjacent matrices in the sample path with the objective of maximizing the correlation between corresponding singular vectors of the two matrices, and (iii) an efficient algorithm is given for untangling the singular value sample paths.

Random matrix processes of the type arising in Multiple Input Multiple Output (MIMO) radio channels will be studied as an example throughout this paper. This is a useful case to study because it lends an intuitive physical interpretation of the SVD and because it is a relevant research area. Spatial multiplexing of information across antennas in an array allows MIMO transceivers to communicate information over independent parallel channels. The SVD of a MIMO channel matrix is used to provide the transmitter and receiver with the necessary beamforming coefficients for spatial multiplexing \cite{Andersen2000}.

The remainder of the work is organized as follows. In Section II, the traditional definition of the SVD is revisited and MIMO beamforming is introduced. Section III details the approach for computing the SVD's of two adjacent matrices and the untangling algorithm is presented. Section IV demonstrates the effectiveness of the algorithm by presenting results from the untangling of a MIMO Gaussian random process. A summary of the major contributions is given in Section V.
 
\section{SVD and Matrix Sample Paths}
In this section, the SVD will be reformulated to allow greater freedom when computing the factors of a random matrix. Next, it will be shown how the SVD of a MIMO channel matrix is used to form a spatial filter that is matched to the radio propagation paths between transmitter and receiver arrays. The problem of tangled sample paths will then be introduced.
\subsection{The SVD Revisited}
Consider the set $\left\{\left.H^{(k)}\in\mathbb{C}^{M\times N}\right|\,k=1,\ldots,K\right\}$. Each $H^{(k)}$ admits a SVD according to,
	\begin{equation}\label{standard_svd}
		H^{(k)}=U^{(k)}\Sigma^{(k)}V^{(k)^H}\text{,}
	\end{equation}
where $U^{(k)}\in\mathbb{C}^{M\times M}$ and $V^{(k)}\in\mathbb{C}^{N\times N}$ are unitary matrices whose columns, $\left\{\vec{u}_1^{\left(k\right)},\ldots, \vec{u}_M^{\left(k\right)}\right\}$ and
$\left\{\vec{v}_1^{\left(k\right)},\ldots, \vec{v}_N^{\left(k\right)}\right\}$, are the left and right singular vectors of $H^{\left(k\right)}$. The matrix, $\Sigma^{(k)}\in\mathbb{R}^{M\times N}$, is populated by the singular values of $H^{\left(k\right)}$, %
$\left\{\sigma_1^{(k)},\ldots,\sigma_{\min\left(M,N\right)}^{(k)}\right\}$, %
along the diagonal and zeros elsewhere. Two \textsl{strict identification conditions} are imposed on the SVD:
	\begin{enumerate}
		\item[i.] The singular values are real and non-negative.
		\item[ii.] The singular values are ordered along the diagonal as,
			\begin{equation}\label{magnitude_sval_ordering}
				\sigma_1^{(k)}\geq\sigma_2^{(k)}\geq\cdots\geq\sigma_{\min\left(M,N%
				\right)}^{(k)}
			\end{equation}
		and the singular vectors are ordered accordingly.
	\end{enumerate}
These strict identification conditions may be relaxed as follows:
	\begin{enumerate}
		\item[(a)] It is permissible to multiply
$\vec{u}_i^{\left(k\right)}$ by $\exp\left(j\theta_{U,i}^{(k)}\right)$ and
$\vec{v}_i^{\left(k\right)}$ by $\exp\left(j\theta_{V,i}^{(k)}\right)$ 
without violating (\ref{standard_svd}) so long as $\sigma_i^{(k)}$ is multiplied by
$\exp\left(-j\left(\theta_{U,i}^{(k)}+\theta_{V,i}^{(k)}\right)\right)$.
		\item[(b)] The singular values may be permuted on the diagonal of
$\Sigma^{\left(k\right)}$ without violating (\ref{standard_svd}) so long as the columns of
$U^{(k)}$ and $V^{(k)}$ are permuted accordingly.
	\end{enumerate}	
In what follows, the case of square matrices ($M=N)$ may be
considered without loss of generality (the generalization to rectangular
matrices is easily accommodated). 
Points (a) and (b) can be incorporated in the SVD definition. Let $P^{\left(k\right)}\in\mathbb{C}^{N\times N}$ be a unitary permutation matrix
such that post multiplication of a matrix by $P^{\left(k\right)}$ causes the
reordering of the columns of the matrix. Let
$\Theta^{(k)}\in\mathbb{C}^{N\times N}$ be a diagonal unitary matrix given by,
\begin{equation}\label{theta_matrix}
\Theta^{(k)}=\text{diag}
\begin{bmatrix}
e^{j\theta_1^{(k)}}&e^{j\theta_2^{(k)}}&\cdots&e^{j\theta_N^{(k)}}
\end{bmatrix}\text{,}
\end{equation}
where $\theta_i^{(k)}$ has its domain on $[0,2\pi)$. By reforming the SVD
factors as,
\begin{equation}\label{generic_svd_factors}
\begin{aligned}[b]
U_{\Theta P}^{(k)}&=U^{(k)}\Theta_U^{\left(k\right)}P^{\left(k\right)},\\
\Sigma_{\Theta
P}^{(k)}&=P^{\left(k\right)^H}\Theta_U^{\left(k\right)^H}\Sigma^{(k)}%
\Theta_V^{\left(k\right)}P^{\left(k\right)},\\
V_{\Theta P}^{(k)}&=V^{(k)}\Theta_V^{\left(k\right)}P^{\left(k\right)},
\end{aligned}
\end{equation}
the SVD becomes,
\begin{equation}\label{generic_svd}
H^{(k)}=U_{\Theta P}^{(k)}\,\,\Sigma_{\Theta P}^{(k)}V_{\Theta P}^{(k)^H}.
\end{equation}
There are therefore an infinite number of valid decompositions for each $H^{(k)}$.
\subsection{SVD Beamforming}
The noise free MIMO system is modeled as,
\begin{equation}\label{mimo_system_model}
\begin{aligned}[t]
\vec{y}\,^{\left(k\right)}&=H^{\left(k\right)}\vec{x}\,^{\left(k\right)}\\
&=\left(U^{\left(k\right)}\Sigma^{\left(k\right)}V^{\left(k\right)^H}\right)\vec{x}\,^{\left(k\right)}
\end{aligned}
\end{equation}
where $\left(\vec{x},\vec{y}\right)\in\mathbb{C}^{N\times 1}$
are the transmitted and received signal vectors respectively. A treatment of the noisy case is beyond the scope of this paper and is not necessary for conveying the ideas presented hereafter. A MIMO channel can be resolved into a set of parallel independent channels if the transmitted symbols are spatially multiplexed and the received symbols are spatially demultiplexed as follows
\cite{Telatar1995},
\begin{equation}\label{mux_demux_mimo_model}
\begin{aligned}[b]
\vec{y}\,^{\left(k\right)}&=U^{\left(k\right)^H}H^{\left(k\right)}\left(V^{\left(k\right)}%
\vec{x}\,^{\left(k\right)}\right)\\
&=\Sigma^{\left(k\right)}\vec{x}\,^{\left(k\right)}
\end{aligned}
.
\end{equation}
Furthermore, spectral efficiency is maximized by waterfilling transmit power across the elements of $\vec{x}\,^{\left(k\right)}$ according the channel's singular values. Such a scheme is accomplished by computing the SVD of $H^{\left(k\right)}$ at the receiver and feeding back $\Sigma^{\left(k\right)}$ and $V^{(k)}$ to the transmitter before performing the waterfilled mux/demux communication of $\vec{x}\,^{\left(k\right)}$. The feedback rate must be fast enough so that this can be accomplished before the channel has changed significantly. On the other hand, feedback requires use of the channel (overhead) which results in a loss of spectral efficiency. The optimal feedback rate is therefore a function of the singular values' temporal coherence and feedback overhead. A similar situation arises in beamforming for MIMO OFDM systems where the singular values of each subcarrier channel vary across frequency and optimal pilot density is a function of coherence bandwidth and pilot overhead.
\begin{figure}
\centering
\includegraphics[width=3.45in]{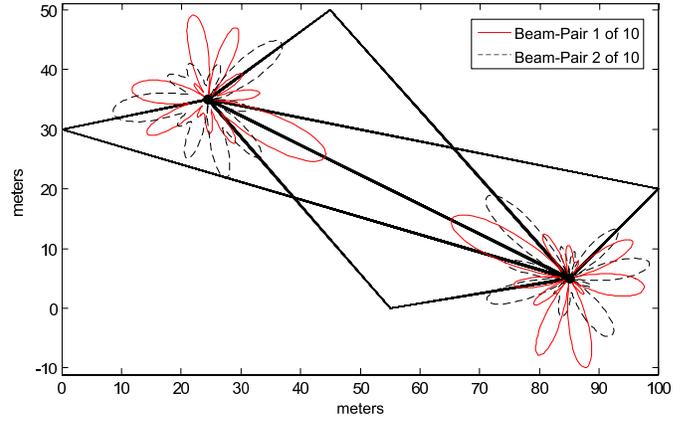}
\caption{Array beam-pairs associated with the two dominant singular channels of a $10\times10$ MIMO channel. Beam  gain-patterns are plotted on a linear scale.}
\label{eigen_beams}
\end{figure}
The mux/demux operations have an important beamforming interpretation \cite{Andersen2000}. Consider the simulated MIMO propagation scenario shown in Fig. (\ref{eigen_beams}). Two circular arrays with ten elements and $1\lambda$ radii are set up in a simple Line of Sight (LoS) scenario with four scatterers. Propagation paths are discovered by raytracing between each pair of transmit and receive antenna. A simple path loss model is then applied to each path and the MIMO channel matrix and its SVD are computed. Pairs of left and right singular vectors are then used as weighting vectors to beamform at the transmitter and receiver \cite{Krim1996}. The beam patterns associated with the two strongest singular values are shown. The first beam-pair focuses energy along the LoS path while the second beam-pair uses the scattering paths. The remaining beam-pair patterns (not shown) are unique and use the available multipath to provide other independent channels. Thus, the SVD relates the MIMO channel matrix to the underlying physical propagation environment by way of the spatial filters formed by the singular vector pairs. It is important to recognize that a common phase rotation of the elements of a singular vector does not alter the associated beam's gain-pattern. Thus, even though there are infinite valid factorizations of a MIMO channel according to (\ref{generic_svd}), there is a unique set of beam gain-patterns for a given propagation scenario.
\subsection{Tangled Sample Paths}
Consider a series, $H^{\left(1\ldots K\right)}=\left\{H^{(1)},\ldots,H^{(K)}\right\}$, that is a sample path taken from a wide sense stationary and bandlimited random matrix process at a rate satisfying the Nyquist-Shannon sampling criteria. Such a sample path arises when sampling a MIMO radio channel across time or frequency. The SVD may be computed for each $H^{\left(k\right)}$ to yield a sample path decomposition,
	\begin{equation}\label{sample_path_decomposition}
	H^{\left(1\ldots K\right)}\longrightarrow U^{\left(1\ldots K\right)},\Sigma^{\left(1\ldots K\right)},V^{\left(1\ldots
	K\right)}
	\end{equation}
where the sample paths of the SVD factors are,
\begin{equation}\label{svd_factor_sample_paths}
\begin{aligned}[b]
&U^{\left(1\ldots K\right)}=\left\{U^{(1)},\ldots,U^{(K)}\right\},\\
&\Sigma^{\left(1\ldots
K\right)}=\left\{\Sigma^{(1)},\ldots,\Sigma^{(K)}\right\},\\
&V^{\left(1\ldots K\right)}=\left\{V^{(1)},\ldots,V^{(K)}\right\},
\end{aligned}
\end{equation}
and sample paths of the singular channels are denoted by the triplet,
\begin{equation}\label{singular_channel_sample_paths}
\left(\vec{u},\sigma,\vec{v}\right)_i^{\left(1\ldots
K\right)}\,\,\,,\,\,\text{for }i=1,\ldots,\min\left(M,N\right).
\end{equation}
A SISO channel sample path, $h_{ij}^{\left(1\ldots K\right)}$, has a smooth evolution that is governed by the covariance of the underlying random process \cite{Jakes1974}. This random process arises from the physical process of multipath superposition. Since both SISO channels and MIMO singular channels arise from similar physical processes, it is expected that the singular channel sample paths should have a smooth evolution. However, when the strict identification conditions, (i) and (ii), are imposed independently from one sample to the next when computing (\ref{sample_path_decomposition}), the resulting singular channel sample paths will be a tangled ordering of the natural singular channel sample paths. This will be made clear by example in subsequent sections. The tangled singular channel sample paths will be seen to not have a smooth evolution and will have an auto-covariance and a cross-covariance that are not consistent with the underlying random process. 
This tangling is a serious problem for closed loop MIMO communications because the feedback rate (pilot density) needed for channel tracking (channel estimation) is greater than what is inherently necessary. Consequently, the system's performance is diminished. 
\section{Untangling Singular Channels}
A solution to the problem of tangled singular channel sample paths can be formulated by using the SVD with relaxed identification conditions to compute singular channel sample paths whose covariance is consistent with that of the underlying random process. The SVD formulation in (\ref{generic_svd}) requires a search over the space of permutation and rotation matrices for the factorization that is most consistent with those of the adjacent matrices' SVDs. Since it is the singular vectors that contain information about the propagation multipath, they will play the central role in finding the optimal factorization. As the search proceeds from one sample to the next, the singular channel sample paths will untangle from the muddled sample paths computed with the strict SVD identification conditions.
\subsection{Correlation Recovery Strategy}
A reference matrix, $H^{\left(\text{R}\right)}$, and a target matrix, $H^{\left(k\right)}$, are chosen from the sample path, $H^{\left(1\ldots K\right)}$. The assignment $H^{\left(\text{R}\right)}=H^{\left(k-1\right)}$ may be made without loss of generality. The SVDs of
$H^{\left(\text{R}\right)}$ and $H^{\left(k\right)}$ are then computed
subject to the strict identification conditions. The triplet
$\left(U^{\left(\text{R}\right)},\Sigma^{\left(\text{R}\right)},V^{\left(
\text{R}\right)}\right)$ will serve as the template to which the triplet
$\left(U^{\left(k\right)},\Sigma^{\left(k\right)},V^{\left(k\right)}\right)$
is matched as closely as possible using permutations and phase rotations according to (\ref{generic_svd_factors}). 
Let the difference matrices between the singular vectors of the
$H^{\left(\text{R}\right)}$ and $H^{\left(k\right)}$ be,
\begin{equation}\label{difference_matrices}
D_U^{\left(k\right)}=U^{\left(\text{R}\right)}-U_{\Theta
P}^{\left(k\right)}\,\,\,,\,\,\,D_V^{\left(k\right)}=V^{\left(\text{R}
\right)}-V_{\Theta P}
\end{equation}
The matching can be formulated as a search over all possible triplets
$\left(\Theta_U^{\left(k\right)},\Theta_V^{\left(k\right)},P^{\left(k\right)}
\right)$ with the objective of minimizing the sum of squared differences,
\begin{equation}\label{min_sum_of_squared_differences}
\min_{P^{\left(k\right)}\text{,
}\Theta_U^{\left(k\right)},\Theta_V^{\left(k\right)}}\,\,\text{tr}\left[D_U^{
\left(k\right)^H}D_U^{\left(k\right)}\right]+\text{tr}\left[D_V^{\left(k
\right)^H}D_V^{\left(k\right)}\right]
\end{equation}
where the sum of squared differences can be expanded as,
\begin{equation}\label{expanded_sum_of_squared_differences}
\begin{aligned}[t]
\text{tr}&\left[D_U^{\left(\text{R}\right)^H}D_U^{\left(k\right)}\right]+
\text{tr}\left[D_U^{\left(\text{R}\right)^H}D_U^{\left(k\right)}\right]\\
&=\text{tr}\left[U^{\left(\text{R}\right)^H}U^{\left(\text{R}\right)}
\right]-2\mathcal{R}e\left\{\text{tr}\left[U^{\left(\text{R}\right)^H}U_{
\Theta P}^{\left(k\right)}\right]\right\}\\
&+\text{tr}\left[U_{\Theta P}^{\left(k\right)^H}U_{\Theta
P}^{\left(k\right)}\right]+\text{tr}\left[V^{\left(\text{R}\right)^H}V^{
\left(\text{R}\right)}\right]\\
&-2\mathcal{R}e\left\{\text{tr}\left[V^{\left(\text{R}\right)^H}V_{\Theta
P}^{\left(k\right)}\right]\right\}+\text{tr}\left[V_{\Theta
P}^{\left(k\right)^H}V_{\Theta P}^{\left(k\right)}\right].
\end{aligned}
\end{equation}
Since terms one, three, four and six on the right side of (\ref{expanded_sum_of_squared_differences})
are not affected by
$\left(\Theta_U^{\left(k\right)},\Theta_V^{\left(k\right)},P^{\left(k\right)}
\right)$, minimizing the sum of squared differences is equivalent to
maximizing the sum of singular vector correlations (terms two and five),
\begin{equation}\label{maximize_correlation}
\max_{P^{\left(k\right)},\Theta_U^{\left(k\right)},\Theta_V^{\left(k\right)}}\,\,\mathcal{R}e\left\{
\text{tr}\left[R_U^{\left(k\right)}+R_V^{\left(k\right)}\right]\right\}
\end{equation}
where the correlation matrices,
\begin{equation}\label{correlation_matrices}
R_U^{\left(k\right)}=U^{\left(\text{R}\right)^H}U_{\Theta
P}\,\,\,,\,\,\,R_V^{\left(k\right)}=V^{\left(\text{R}\right)^H}V_{\Theta P}
\text{,}
\end{equation}
have elements that comprise the set of all possible inner products of the reference and target singular vectors.
The maximization in (\ref{maximize_correlation}) implies a search over the infinite set
$\left\{\left(\Theta_U^{\left(k\right)},\Theta_V^{\left(k\right)},P^{\left(k\right)}\right)\right\}$. Fortunately, the global optimum solution,
$\Big(\widehat{\Theta}\,\!_U^{\left(k\right)},\widehat{\Theta}\,\!_V^{\left(k\right)},\widehat{P}\,\!^{\left(k\right)}\Big)$, can be found by searching over a finite set if either
$\widehat{P}\,\!^{\left(k\right)}$ or
$\Big(\widehat{\Theta}\,\!_U^{\left(k\right)},\widehat{\Theta}\,\!_V^{\left(k\right)}\Big)$ are known. As will be shown, an elegant search strategy is to find
$\widehat{P}\,\!^{\left(k\right)}$ assuming that an appropriate
$\Big(\widehat{\Theta}\,\!_U^{\left(k\right)},\widehat{\Theta}\,\!_V^{\left(k\right)}\Big)$ exists and then to compute
$\Big(\widehat{\Theta}\,\!_U^{\left(k\right)},\widehat{\Theta}\,\!_V^{\left(k\right)}\Big)$ analytically given
$\widehat{P}\,\!^{\left(k\right)}$.

The search for $\widehat{P}\,\!^{\left(k\right)}$ proceeds as follows. A simple interchange columns $i$ and $j$ of $U_{\Theta P}^{\left(k\right)}$ is equivalent to interchanging columns $i$ and $j$ of
$R_U^{\left(k\right)}$ (and similarly for $V_{\Theta P}^{\left(k\right)}$ and
$R_V^{\left(k\right)}$). Consequently, the measure of improvement in, 
\begin{equation}\label{trace_of_correlation_matrix}
\mathcal{R}e\left\{\text{tr}\left[R_U^{\left(k\right)}+R_V^{\left(k\right)}
\right]\right\}=\mathcal{R}e\left\{\sum_{i=1}^n\left(r_{ii}^{\left(k\right)}
\right)_U+\left(r_{ii}^{\left(k\right)}\right)_V\right\}\text{,}
\end{equation}
that results from an interchange of columns $i$ and $j$ is, 
\begin{equation}\label{measure_of_improvement}
\begin{aligned}[b]
\Delta_R^{\left(k\right)}\left(i,j\right)=&\mathcal{R}e\Bigl\{\left(r_{ij}^{
\left(k\right)}+r_{ji}^{\left(k\right)}-r_{jj}^{\left(k\right)}-r_{ii}^{
\left(k\right)}\right)_U\\
&+\left(r_{ij}^{\left(k\right)}+r_{ji}^{\left(k\right)}-r_{jj}^{\left(k
\right)}-r_{ii}^{\left(k\right)}\right)_V\Bigr\}
\end{aligned}.
\end{equation}
Phase rotations of the singular vectors according to (\ref{generic_svd_factors}) by the appropriate
$\left(\Theta_U^{\left(k\right)},\Theta_V^{\left(k\right)}\right)$ guarantee that the greatest possible improvement,
\begin{equation}\label{greatest_possible_improvement}
\begin{aligned}[b]
\Delta_{\left|R\right|}^{\left(k\right)}\left(i,j\right)&=\left(
\left|r_{ij}^{\left(k\right)}\right|+\left|r_{ji}^{\left(k\right)}\right|-
\left|r_{jj}^{\left(k\right)}\right|-\left|r_{ii}^{\left(k\right)}\right|
\right)_U\\
&+\left(\left|r_{ij}^{\left(k\right)}\right|+\left|r_{ji}^{\left(k\right)}
\right|-\left|r_{jj}^{\left(k\right)}\right|-\left|r_{ii}^{\left(k\right)}
\right|\right)_V
\end{aligned}
\text{,}
\end{equation}
is achievable, thus permitting the search for
$\widehat{P}\,\!^{\left(k\right)}$ to be executed by column swaps
without jointly searching for
$\widehat{\Theta}\,\!_U^{\left(k\right)}$ and
$\widehat{\Theta}\,\!_V^{\left(k\right)}$. The maximization problem
is then,
\begin{equation}\label{simplified_maximization}
\max_{P^{\left(k\right)}}\,\,\text{tr}\left[\left|R_U^{\left(k\right)}
\right|+\left|R_V^{\left(k\right)}\right|\right].
\end{equation}
The search for $\widehat{P}\,\!^{\left(k\right)}$ can be performed
iteratively by column swaps using (\ref{greatest_possible_improvement}) as a metric for directing the search.
The details of an efficient algorithm for implementing this search strategy
are given in Section III-B.
Once $\widehat{P}\,\!^{\left(k\right)}$ has been found, finding the
associated
$\Big(\widehat{\Theta}\,\!_U^{\left(k\right)},\widehat{\Theta}\,\!_V^{\left(k\right)}\Big)$ is simply a matter of computing the angle between singular vectors of $H^{\left(\text{R}\right)}$ and the permuted singular vectors of $H^{\left(k\right)}$ according to,
\begin{equation}\label{angle_btwn}
\left\{\left.\theta_{U,i}^{\left(k\right)}=-\arg\left(\vec{u}_i^{\left(
\text{R}\right)^H}\vec{u}_i^{\left(k\right)}\right)\right|i=1,\ldots,N\right
\}\text{,}
\end{equation}
and likewise for $\big\{\theta_{V,i}^{\left(k\right)}\big\}$ where $\arg\left(\cdot\right)\equiv\mathcal{I}m\left\{\ln\left(\cdot\right)\right\}$ yields the radian angle of a complex scalar. These left and right phase rotation matrices are then applied to $\Sigma^{(k)}$ as required by (\ref{generic_svd}) to compensate for the conjugate phase rotation being applied to the corresponding singular vectors. The resulting match of $\big(U_{\Theta P}^{\left(k\right)}\,,V_{\Theta P}^{\left(k\right)}\big)$ to $\left(U^{\left(\text{R}\right)}\,,V^{\left(\text{R}\right)}\right)$ is
optimal in a least squares sense.

Thus far, the discussion has assumed noise-free sample paths. The case of an additive noise model, %
$\vec{y}\,^{\left(k\right)}=H^{\left(k\right)}\vec{x}\,^{\left(k\right)}+\vec{n}\,^{\left(k\right)}$%
, is relevant and can be addressed in part by using the weighted singular vectors,
\begin{equation}
\begin{aligned}[t]
&\widetilde{U}^{\left(k\right)}=U^{\left(k\right)}\left(\Sigma^{%
\left(k\right)}\right)^{1/2}\,,\,\widetilde{V}^{\left(k%
\right)}=V^{\left(k\right)}\left(\Sigma^{\left(k\right)}\right)^{1/2}
\end{aligned}
\end{equation}
in the search for
$\left(\widehat{\Theta}_U^{\left(k\right)},%
\widehat{\Theta}_V^{\left(k\right)},%
\widehat{P}^{\left(k\right)}\right)$. This biases the search in favor of those singular vectors
with greatest signal to noise ratio. This is not all that can be done to
counteract the noise. However, a proper treatment of noisy sample paths is
beyond the scope of this paper.
\subsection{Untangling Algorithm}
The algorithm given in Table 1 untangles the singular channel sample paths from the sample paths computed with the strict identification conditions. The algorithm's core (lines 4-21) is an efficient implementation of the correlation recovery strategy discussed in Sections III-A. The number of
possible permutations on the singular vector ordering is $N!$. However, the core algorithm guarantees that the least squares solution for $\widehat{P}^{\left(k\right)}$ is found in at most $N\left(N-1\right)$ column swaps. 
The essential aspects of the core algorithm are: 
\begin{figure}
\centering
\includegraphics[width=3.0in]{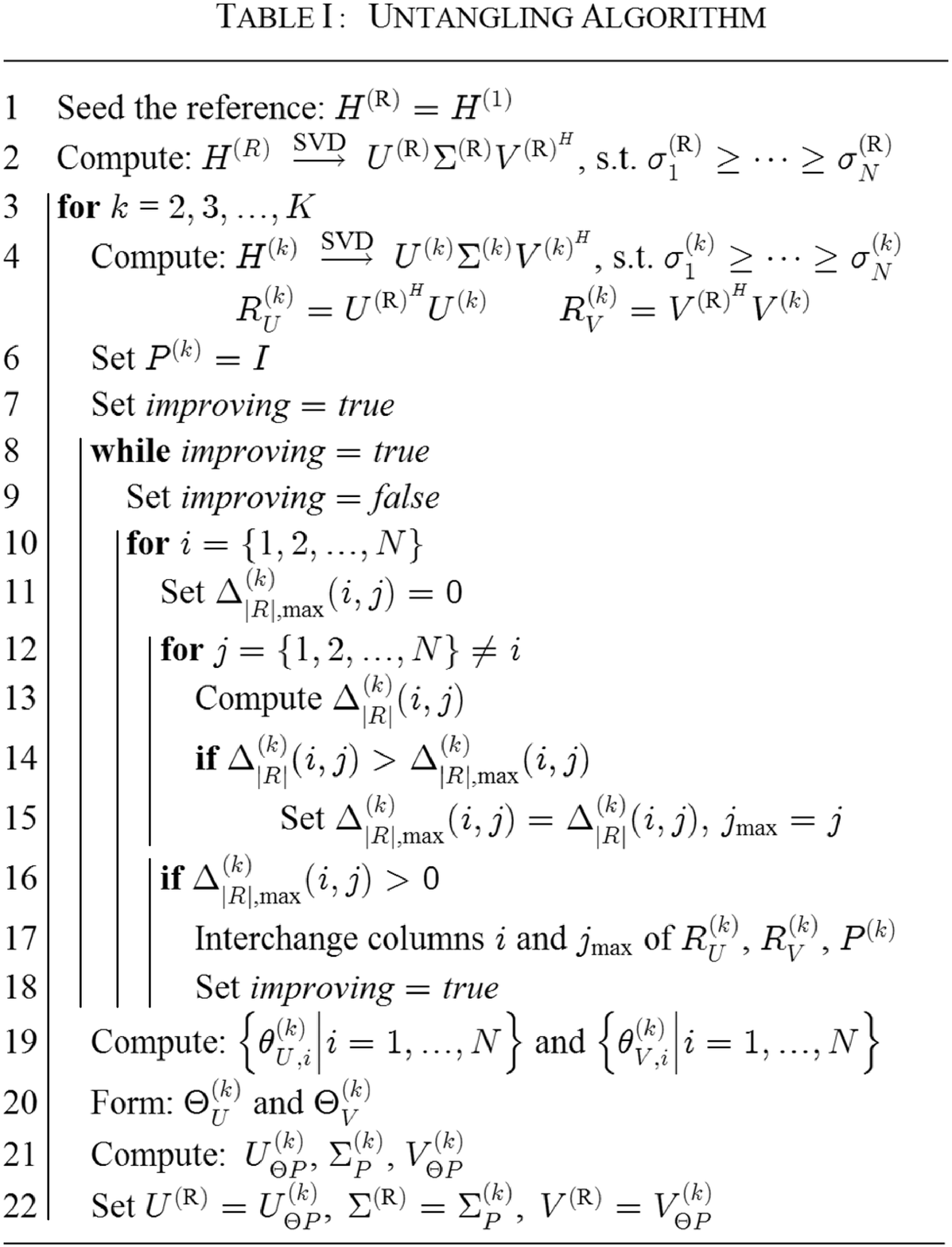}
\label{pseudocode}
\end{figure}
	\begin{enumerate}
		\item[i.] Column $i$ (pivot column) is swapped with the column that would
result in the maximum improvement to (\ref{greatest_possible_improvement}), whenever an improvement is
possible.
		\item[ii.] Each column becomes the pivot column in cyclic order and the
algorithm proceeds iteratively until no further improvement can be made.
	\end{enumerate}
More generally, the core algorithm can be applied to any problem where the trace of a matrix must be maximized under column (or row) permutations.
The algorithm's shell (lines 1-3,22) parses the sample path, $H^{\left(1 \ldots K \right)}$, using a sliding reference, $H^{\left(\text{R}\right)}=H^{\left(k-1\right)}$. The SVD factors computed by the core algorithm for each matrix sample are therefore dependent on the reference seed-matrix's factorization (line 2). This seed matrix has its SVD computed according to strict identification conditions on ordering and phase. Any ordering and phase for the seed reference's SVD would also suffice. Whatever the initial order and phase, they determine the relative ordering and phase-offset of the untangled singular channel sample paths. As will be discussed in Section IV, the untangling solution is unique when conditioned on the seed matrix's singular value ordering and phases.

The Untangling Algorithm's performance can be poor if either of the following two situations arise and are not treated appropriately. The first situation is that of a matrix with close singular values. In this case, the SVD algorithm becomes unstable and slight perturbations of the matrix cause wildly different decompositions. The second situation (an extreme case of the first) is that of a matrix with repeated singular values. In this case, the singular vectors associated with the repeated singular values are indeterminate and contain no information. The case of equal singular values is very unlikely in practical applications where the $H^{\left(k\right)}$ are sampled from a random processes occurring in the physical world. These two situations may be handled by implementing the following simple modifications to the Untangling Algorithm:
	\begin{enumerate}
		\item[i.] Set a threshold, $\mu^{\left(k\right)}$, for the tolerable difference between singular values returned by the SVD algorithm. For example,
			\begin{equation}\label{threshold}
				\min\left(\sigma_i^{\left(k\right)}-\sigma_{i+1}^{\left(k\right)}\right)<\epsilon,i=1,\ldots,N-1
			\end{equation}
		where $\epsilon$ is a small positive number and %
		$\sigma_1^{\left(k\right)}\geq\cdots\geq\sigma_N^{\left(k\right)}$.
		\item[ii.] If the threshold is violated then attempt a match but display a
warning message and blacklist the offending matrix so that it cannot be used
as a reference.
	\end{enumerate}
It is possible, though unlikely, that a tie is computed for different column swaps' improvement metrics. This situation is averted in line 14 by the greater-than logic which causes the first of the column swaps ties to be used. Even so, a tie should be treated in a manner similar to the case of repeated singular values by giving a warning, attempting a match and blacklisting the offending matrix. 
\section{Results}
In this section, an example is given to demonstrate the effectiveness of the Untangling Algorithm in computing smoothly evolving singular channels from a fading MIMO radio channel. A comparison is made between tangled and untangled sample paths. The covariance of the untangled singular value sample path is compared with those of the SISO channel sample paths and tangled singular value sample paths. Joint and marginal distributions of the untangled singular values are compared with those known from random matrix theory.
\subsection{MIMO Channel Synthesis}
Consider the case of a MIMO radio distortion process between a fixed transmitter and a mobile receiver. It may be assumed that the random MIMO distortion process, $\bs{H}$, is a time varying matrix of i.i.d. zero mean Gaussian random variables, $\bs{h}_{ij}$. The temporal auto-covariance of each SISO channel realized from $\bs{H}$ may be given by \cite{Jakes1974},
\begin{equation}\label{bessel_correlation}
c\left(h_{ij}^{\left(t\right)},h_{ij}^{\left(t+\tau\right)}\right)=J_0\left(2\pi f_d\,\tau\right)
\end{equation}
where $J_0\left(x\right)$ is the Bessel function of the first kind and $f_d=v/\lambda_c$ is the maximum Doppler frequency resulting from a receiver moving with velocity, $v$, and a carrier frequency with wavelength, $\lambda_c$. The assumption that the SISO channels are independent and identically distributed is justifiable in many communication scenarios. The temporal auto-covariance given in (\ref{bessel_correlation}) arises from non-line-of-sight scenarios with rich multipath propagation. This bandlimited random matrix process may be sampled over time to yield a sample path, $H^{\left(1\ldots K\right)}$. Realizations of $h_{ij}^{\left(1\ldots K\right)}$ can readily be generated in computer simulation by a multipath sum-of-sinusoids model \cite{Clarke1968}. This model does not generate MIMO channels directly from antenna array geometries and a scattering topology as in the case of Fig. \ref{eigen_beams}. However, array geometries and a scattering topology are associated with the model's output by way of the the random variables used to seed the paths' phase and arrival angle.  For the analysis that follows, MIMO channel sample paths were generated according to the sum-of-sinusoids model with these parameters: $500$ paths, $f_d=15\,$Hz, $\lambda_c=0.125\,$m, and $T_s=3.3\,$ms.
\subsection{Untangled Sample Paths}
\begin{figure}
\centering
\includegraphics[width=3.4in]{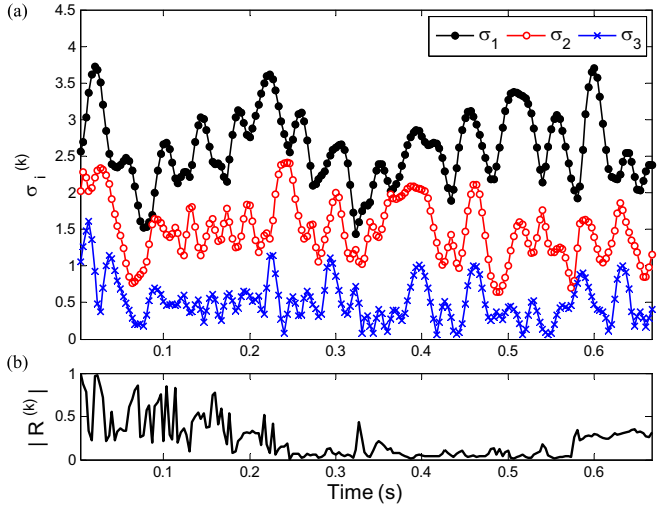}
\caption{(a) Singular value sample paths of a spatially white $3 \times 3$ MIMO channel computed using strict identification conditions on singular value ordering and phase. (b) Average singular vector correlation between adjacent samples.}
\label{iid_raw_svals}
\end{figure}

\begin{figure}
\centering
\includegraphics[width=3.4in]{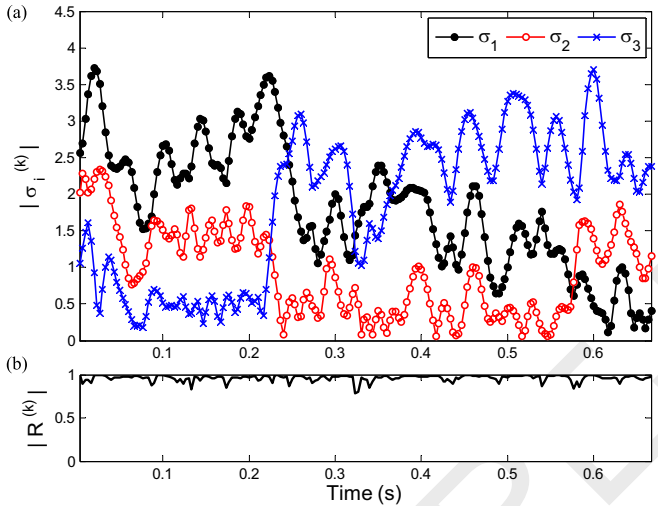}
\caption{(a) Magnitude of the untangled singular value sample paths of a spatially white $3 \times 3$ MIMO channel. (b) Average singular vector correlation between adjacent samples.}
\label{iid_untangled_svals}
\end{figure}

\begin{figure}
\centering
\includegraphics[width=3.4in]{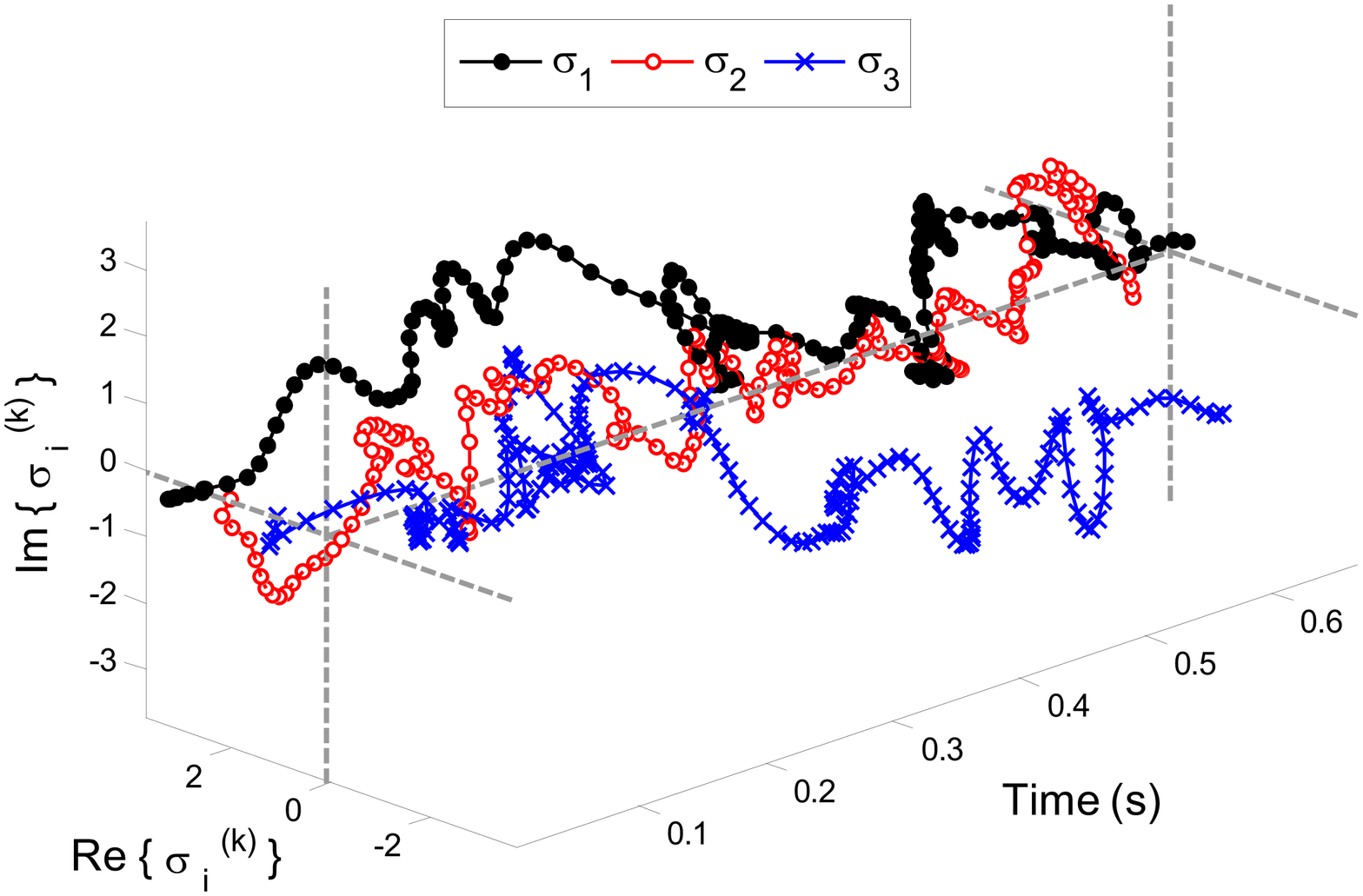}
\caption{Smooth evolution of the untangled singular value sample paths from a spatially white $3 \times 3$ MIMO channel.}
\label{iid_untangled_complex_svals}
\end{figure}

Fig. \ref{iid_raw_svals}(a) shows the singular sample paths for a realization of $H^{\left(1\ldots K\right)}$ where $M=N=3$. SVD's were computed according to strict identification conditions on ordering and phase with the result that singular value sample paths form separate layers. Fig. \ref{iid_raw_svals}(b) shows the average singular vector correlation between adjacent samples,
\begin{equation}\label{svec_correlation}
\overline{R}^{\left(k\right)}=\frac{1}{2N}\text{tr}\left(U^{\left(\text{R}%
\right)^H}U_{\Theta P}^{\left(k\right)}+V^{\left(\text{R}\right)^H}V_{\Theta
P}^{\left(k\right)}\right).
\end{equation}
The correlation profile indicates that although the singular vectors of the $k^{\text{th}}$ and $\left(k-1\right)^{\text{th}}$ samples are initially correlated, the strict identification conditions soon cause a loss of correlation.
Fig. \ref{iid_untangled_svals}(a) shows the \emph{magnitude} of the singular values after applying the Untangling Algorithm to the same MIMO sample path as was used for Fig. \ref{iid_raw_svals}(a). Although the singular values magnitudes have not changed, the traces that connect them are very different than what they were before untangling. The singular value sample paths are now observed to weave. When the traces cross, they do so in natural smooth transitions. The correlation profile in Fig. \ref{iid_untangled_svals}(b) shows that the untangled singular vectors maintain high correlation with those of the earlier sample.
An important aspect of the untangling process is that the singular values are allowed to become complex valued. Fig. \ref{iid_untangled_complex_svals} shows the untangled singular value sample paths as a function of time. These traces reflect a natural smooth evolution in accordance with the bandlimited underlying random process.

Furthermore, the untangling solution is unique when conditioned on the seed matrix's singular value ordering and phases. For example, consider the forward untangled singular value sample path %
$\Sigma_\rightarrow^{\left(1,2,\ldots K\right)}$ %
generated by parsing %
$H^{\left(1,2,\ldots K\right)}$ %
and the reverse untangled singular value sample path, %
$\Sigma_\leftarrow^{\left(1,2,\ldots K\right)}$ %
generated by parsing %
$H^{\left(K\ldots ,2,1\right)}$.
It was verified that,
\begin{equation}\label{forward-reverse}
\Sigma_\rightarrow^{\left(1,2,\ldots K\right)}=
P_r\odot\left(
\Sigma_\leftarrow^{\left(K,\ldots,2,1\right)}
\odot\Theta_r
\right)\odot P_r
\end{equation}
without exception (except for numerical precision errors). Here $\Theta_r$ is a diagonal unitary matrix that derotates the singular values in $\Sigma_\leftarrow^{(K)}$ so they are real and non-negative, $P_r$ is the permutation matrix that is then used to perform the reverse one-to-one mapping of the singular values in $\Sigma_\rightarrow^{(K)}$ to those in $\Sigma_\leftarrow^{(1)}$ and $\odot$ represents element-wise matrix multiplication over the entire sample path. The existence of a unique solution for a given MIMO sample path is consistent with the specific underlying physical process giving rise to that sample path.
\subsection{Sample Path Covariance}
\begin{figure}
\centering
\includegraphics[width=3.4in]{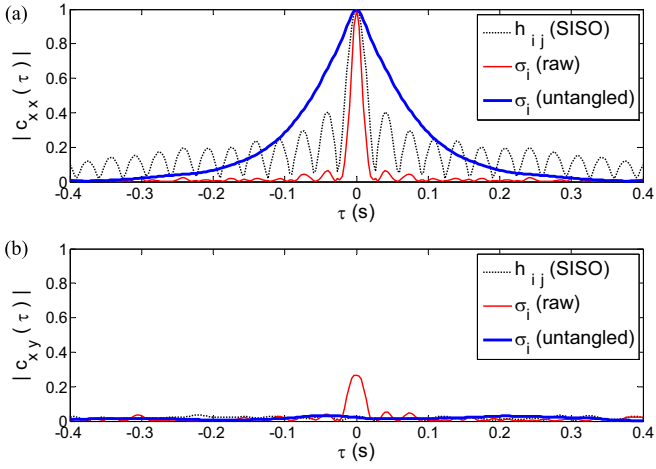}
\caption{Covariance functions computed from a spatially white $3 \times 3$ MIMO channel (a) auto-covariance (b) cross-covariance.}
\label{covariance}
\end{figure}

Fig. \ref{covariance} compares the magnitude of the temporal covariance of a SISO channel sample path, a raw singular value sample path, and an untangled singular value sample path. These are computed from $3 \times 3$ sample paths of length $K=10\,000$. In general, each of the $\min\left(M,N\right)$ raw singular value sample paths has a slightly different covariance function whereas the covariance of each of the untangled singular value sample path is the same. For simplicity, the covariance shown for the raw singular values is the mean of all covariances. 

Figure \ref{covariance}(a) shows that the SISO channel auto-covariance matches that expected from (\ref{bessel_correlation}) and reaches $0.7$ at $\tau=0.012\,$s. The auto-covariance for the raw singular value sample path decays slightly faster than that of the raw singular channels and reaches $0.7$ at $\tau=0.006\,$s.
It is remarkable to observe that the untangled singular channel's covariance decays far less rapidly than that of the SISO channels and reaches $0.7$ at $\tau=0.035\,$s. The untangling process has `slowed down' the perceived channel dramatically. This result implies that channel tracking (estimation) can be done with a lower feedback rate (pilot density) when using the untangled singular channels instead of the raw MIMO channels directly. 
This promise of improved efficiency can be realized by implementing the Untangling Algorithm at the receiver. Furthermore, a transmitter with a priori knowledge of the singular channels' covariance matrix function has an extended horizon for predicting the channel. This is a benefit because latency due to channel estimation and feedback carries the consequence that the transmit beamforming coefficients may no longer be matched to the current channel. The transmitter may then use the most recent untangled channel state information to predict the current channel state for beamforming and waterfilling.

Figure \ref{covariance}(b) compares the magnitude of the temporal cross-covariance of a SISO channel sample path, a raw singular value sample path, and an untangled singular value sample path. As expected, the SISO channel sample paths are uncorrelated. The tangled nature of the raw singular value sample paths manifests itself as a non-zero covariance at lag $\tau=0$. The fact that the untangled singular value sample paths have no covariance verifies that the untangling process has resolved sample paths whose evolutions are independent of each other.
\subsection{Sample Path Density Functions}
Analytic results are known for the distributions of ordered and unordered eigenvalues of the Wishart matrices \cite{Telatar1995}, \cite{Tulino2004}. This is relevant to the topic of MIMO singular channels because $W^{\left(k\right)}=H^{\left(k\right)}H^{\left(k\right)^H}$ is drawn from the Wishart distributed random process $\bs{W}=\bs{HH}^H$. The eigenvalues of $W^{\left(k\right)}$, $\left(\lambda_1^{(k)},\ldots,\lambda_M^{(k)}\right)$, are simply the magnitude squared of the singular values of $H^{\left(k\right)}$. Given that $H^{\left(k\right)}\in\mathbb{C}^{M\times N}$, the joint density function for the unordered eigenvalues of $\bs{W}$ is,
\begin{equation}\label{unordered_jpdf}
\begin{aligned}[b]
P\left(\lambda_1,\ldots,\lambda_M\right)=\frac{1}{M!}\exp\left(-\Sigma_{i=1}^M\lambda_i\right)&%
\prod_{i=1}^M\frac{\lambda_i^{N-M}}{\left(M-i\right)!\left(N-i\right)!}\\&%
\prod_{i<j}^M\left(\lambda_i-\lambda_j\right)^2
\end{aligned}
\text{,}
\end{equation}
while the marginal density function of the unordered eigenvalues of $\bs{W}$ is,
\begin{equation}\label{unordered_mpdf}
P\left(\lambda\right)=\frac{1}{M}\sum_{i=0}^{M-1}\frac{i!}{\left(i+N-M%
\right)!}\left[L_i^{N-M}\left(\lambda\right)\right]^2\frac{\lambda^{M-N}}{e^%
\lambda}
\end{equation}
where $L_n^k\left(x\right)$ is the associated Laguerre polynomial with Rodrigues representation,
\begin{equation}\label{laguerre_polynomial}
L_n^k\left(x\right)=\sum_{m=0}^n\left(-1\right)^m\frac{\left(n+k\right)!}{%
\left(n-m\right)!\left(k+m\right)!\,m!}x^m
\end{equation}

\begin{figure}
\centering
\includegraphics[width=3.4in]{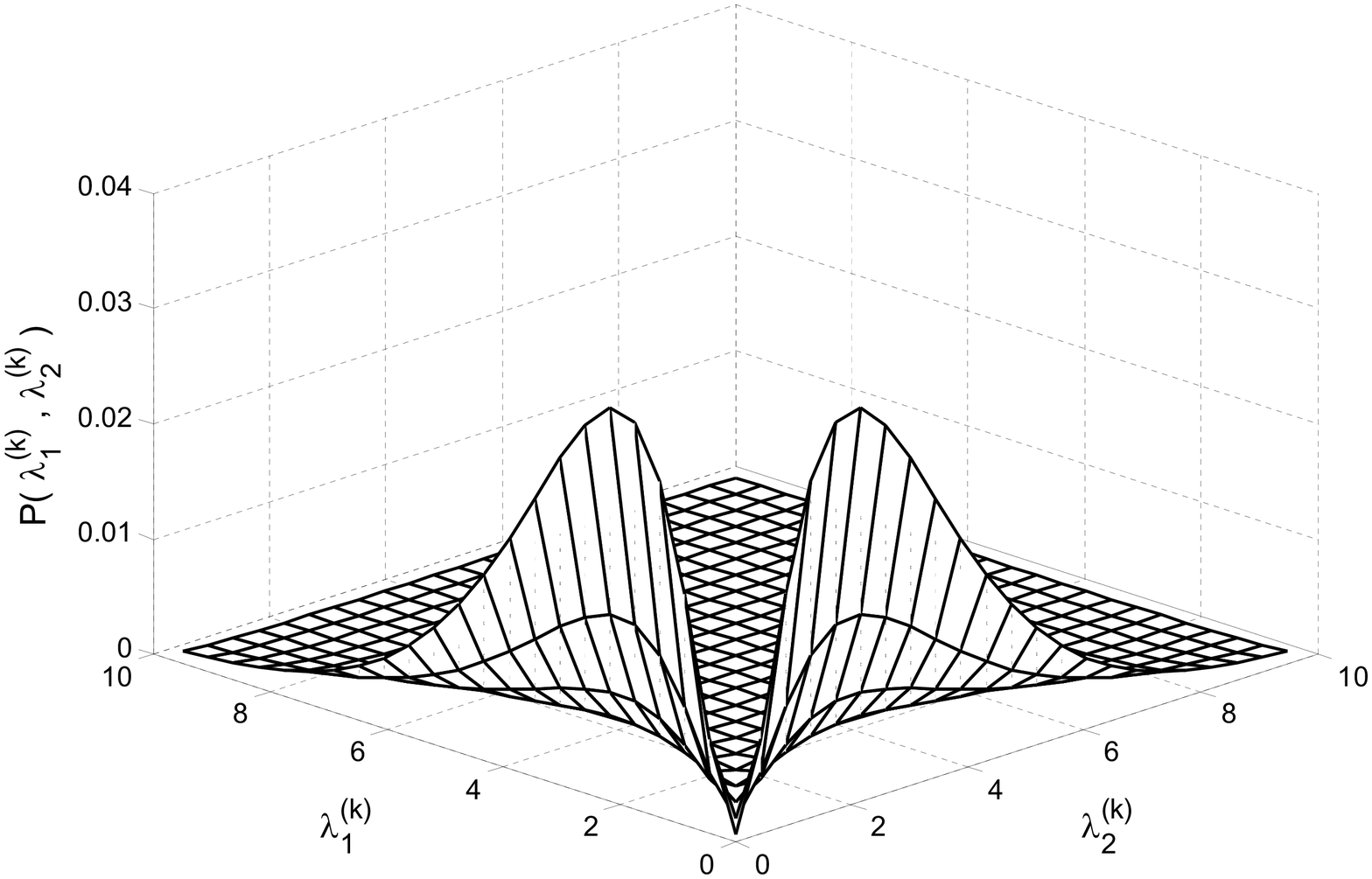}
\caption{Joint density function given by (\ref{unordered_jpdf}) for the unordered eigenvalues of $\bs{W}$ where $\bs{H}\in\mathbb{C}^{2\times2}$ and $\bs{H}\sim\mathcal{N}\left(0,I\right)$}
\label{joint_pdf_theory}
\end{figure}

\begin{figure}
\centering
\includegraphics[width=3.4in]{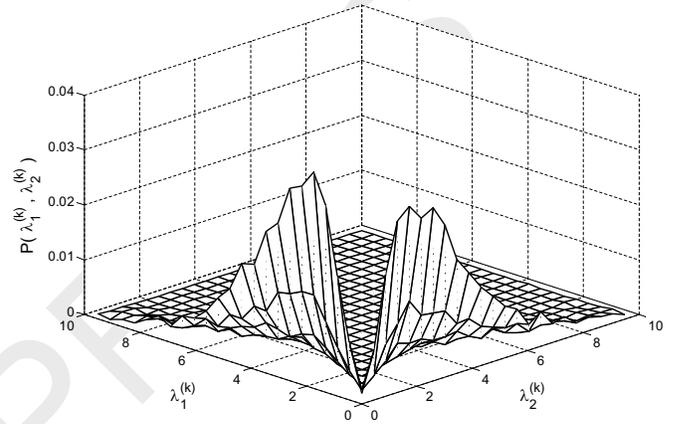}
\caption{Joint distribution of the eigenvalues of $W^{\left(k\right)}$ from an untangled sample path $H^{\left(1\ldots K\right)}$ where $H^{\left(k\right)}\in\mathbb{C}^{2\times 2}$ and $H^{\left(k\right)}\sim\mathcal{N}\left(0,I\right)$.}
\label{joint_pdf_untangled}
\end{figure}

\begin{figure}
\centering
\includegraphics[width=3.4in]{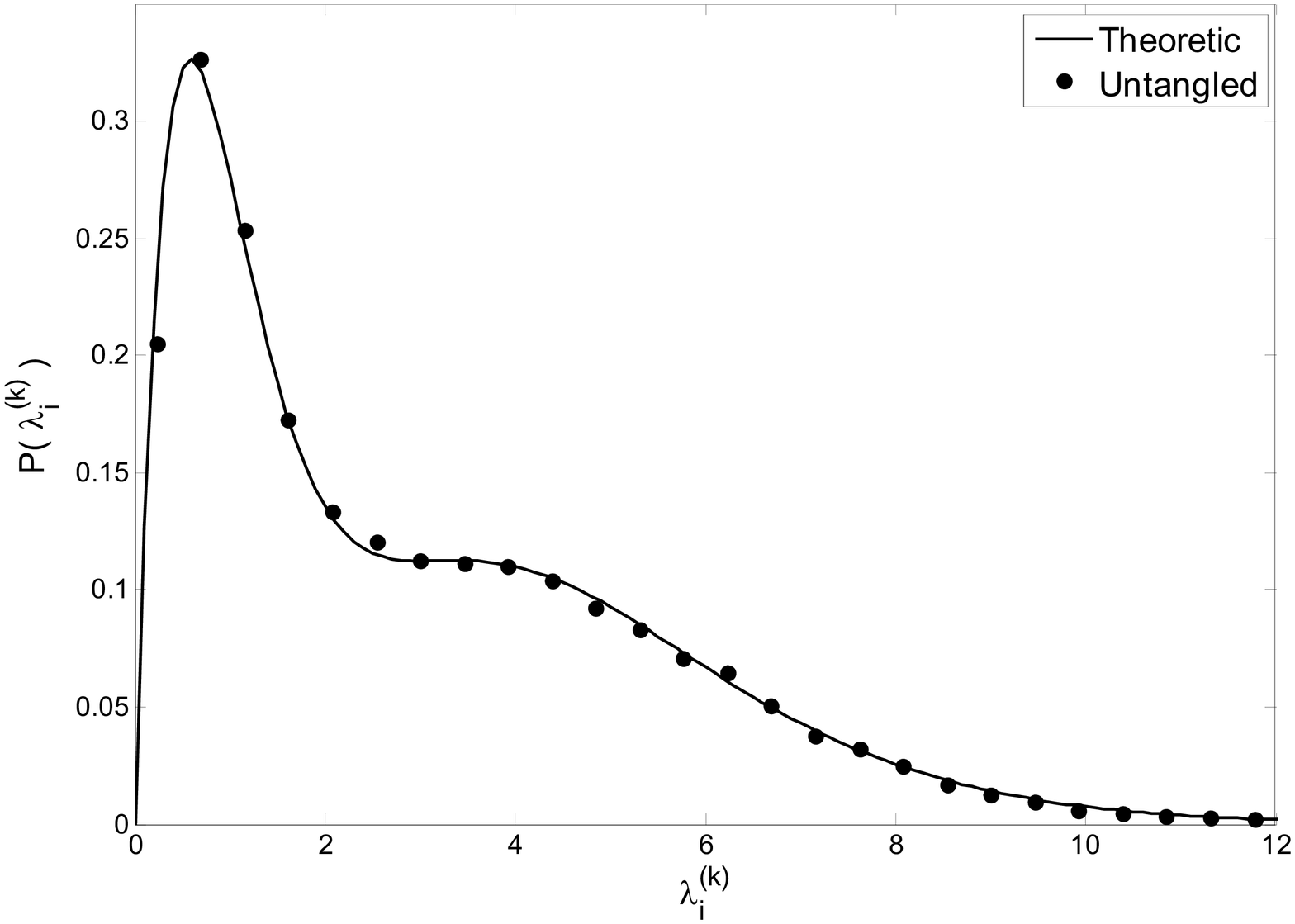}
\caption{(i) Marginal density function given by (\ref{unordered_mpdf}) for unordered eigenvalues of $\bs{W}$ where $\bs{H}\in\mathbb{C}^{2\times 3}$ and $\bs{H}\sim\mathcal{N}\left(0,I\right)$ and (ii) computed eigenvalue marginal distribution for $W^{\left(k\right)}$ from an untangled sample path $H^{\left(1\ldots K\right)}$ where $H^{\left(k\right)}\in\mathbb{C}^{2\times 3}$ and $H^{\left(k\right)}\sim\mathcal{N}\left(0,I\right)$.}
\label{marginal_pdf}
\end{figure}

Fig. \ref{joint_pdf_theory} shows the joint density function computed from (\ref{unordered_jpdf}) for $N=M=2$. Fig. \ref{joint_pdf_untangled} shows the joint distribution of the untangled eigenvalues for a realization of $H^{\left(1\ldots 10000\right)}$ and computed using a histogram over the same domain as the function in Fig. \ref{joint_pdf_theory}. The untangled distribution is matches the theoretic density function and is completely different from the joint distribution of the ordered eigenvalues (not shown). It is difficult to show the joint density function for $\min\left(M,N\right)>2$ but the marginal density function for these higher order cases can be shown as is done in Fig. \ref{marginal_pdf} for $M=2$ and $N=3$. Again, there is very good agreement between the untangled eigenvalue distribution and the theoretic distribution. It is expected that the untangled eigenvalues should have distributions agreeing with those given in (\ref{unordered_jpdf}) and (\ref{unordered_mpdf}) because the untangling operation recovers the  ordering randomness of the underlying i.i.d. Gaussian process.
\section{Conclusion}
This work has addressed the problem of how to compute the SVD's of a random matrix sample path in a way that preserves the covariance of the underlying random process. The solution to this problem was developed as follows: (i) a SVD with relaxed identification conditions was proposed, (ii) an approach was formulated for computing the SVD's of two adjacent matrices in the sample path that maximizes the correlation between corresponding singular vectors of the two matrices, and (iii) an efficient algorithm was given for untangling the singular value sample paths. 

The algorithm's effectiveness was demonstrated on i.i.d. Gaussian MIMO channels. It was shown that the algorithm resolves smoothly evolving singular channel sample paths that are in accord with the stochastic structure of the underlying random process. Furthermore, the algorithm gives a unique solution conditioned on the seed matrix's singular value ordering and phases.

It was shown that the untangling process dramatically increases the coherence period (or bandwidth) of the singular channels. A primary application of the algorithm is in MIMO radio systems. The benefit promised by using SVD untangling in these systems is that the fading rate of the channel's SVD factors is greatly reduced so that the performance of channel estimation, channel feedback and channel prediction can be increased.


\bibliographystyle{IEEEtran}
\bibliography{C:/kingston/(DOCS)_Library/decibel}

\end{document}